\documentclass{article}
\usepackage{graphicx} 

\usepackage[utf8]{inputenc}
\usepackage{ragged2e}
\usepackage{amsfonts}
\usepackage{mathtools}
\usepackage{amssymb}
\usepackage{amsthm}
\usepackage{amsmath}
\usepackage{afterpage}
\usepackage{authblk}

\usepackage{algorithm}
\usepackage{algorithmic} 
\usepackage{color}
\usepackage{hyperref}
\hypersetup{
    colorlinks=true, 
    linktoc=all,     
    linkcolor = blue,
}

\usepackage{geometry}
\newcommand{\mm}{{\overline{m}}}

\newcommand{\PP}{{\overline{P}}}
\newcommand{\QQ}{{\overline{Q}}}
\newcommand{\qq}{{\overline{q}}}
\newcommand{\ZZ}{{\mathbb{Z}}}

\theoremstyle{definition}
\newtheorem{definition}{\bf {Definition}}[section]

\newtheorem{problem}[definition]{\bf {Problem}}
\newtheorem{example}[definition]{\bf {Example}}

\newtheorem{lemma}[definition]{\bf {Lemma}}
\newtheorem{theorem}[definition]{\bf {Theorem}}

\newtheorem{proposition}[definition]{\bf Proposition}
\newtheorem{remark}[definition]{\bf Remark}

\title{Root Extraction in Finite Abelian Groups}
\author{Udvas Acharjee}
\author{M S Srinath}
\affil{Department of Mathematics and Computer Science, Sri Sathya Sai Institute of Higher Learning, Prasanthi Nilayam, Puttaparthi, 515134, Andhra Pradesh, India}

\date{}

\begin{document}

\maketitle
\begin{abstract}
    We formulate the \emph{Root Extraction problem} in finite Abelian $p$-groups and then extend it to generic finite Abelian groups. We provide algorithms to solve them. We also give the bounds on the number of group operations required for these algorithms. We observe that once a basis is computed and the discrete logarithm relative to the basis is solved, root extraction takes relatively fewer "bookkeeping" steps. Thus, we conclude that root extraction in finite Abelian groups is \emph{no harder} than solving discrete logarithms and computing basis. 
\end{abstract}
\section{Introduction}

A simple form of the \emph{root extraction} is as follows:
Let $G  = \langle P,Q \rangle$ where 
\begin{equation}\label{struct1}
    G \approx \frac{\ZZ}{\ell^e\ZZ}\times\frac{\ZZ}{\ell^e\ZZ}
\end{equation}
 Consider an element $K\in G$ and $m,n\in \frac{\ZZ}{\ell^e\ZZ}$ such that,\begin{equation}\label{eqn1}
    K=mP+nQ
\end{equation}

If element $K$ and the multipliers $m,n$ are known, the root extraction problem is to find a basis $P,Q$ of $G$ such that \ref{eqn1} holds. The following table summarizes this problem and related problems.
\begin{table}[htbp]
	\centering
		\begin{tabular}{|l|c|c|c|}
			\hline
			\hline
			Problem & Element & Multipliers & Base Pts. \\
			        &  $K$    &   $m,n$ & $P,Q$   \\
			\hline
			\hline
			Exponentiation & ? & \checkmark & \checkmark \\
			\hline
			Extended DLP & \checkmark & ? & \checkmark \\
			\hline
			Root Extraction & \checkmark & \checkmark & ?\\
			\hline
            Basis computation & - & - & ?\\
			\hline
            \hline
		\end{tabular}
\end{table}

The root extraction problem has been solved for the groups of the form (\ref{struct1})  by Srinath in \cite{MSS}. Similarly for the algorithms to solve the discrete logarithm problem see \cite{Teske} and \cite{Sutherland2008StructureCA}, for the basis computation problem see \cite{Sutherland2008StructureCA}, and the square and multiply algorithm for exponentiation is given in \cite[Chapter 9, \S 9.2]{UdiManberAlgoBook}.
Lower bound for the root extraction problem in generic finite Abelian groups is given in \cite{cryptoeprint:2002/013}.

The method to solve the root extraction problem discussed in \cite{MSS} can be easily extended to any group $G$ of the form:\begin{equation}\label{Struct2}
  G \approx
\underbrace{\frac{\ZZ}{p^{e}\ZZ}\times...\times\frac{\ZZ}{p^e\ZZ}}_{N\text{-times}}
\end{equation}

Our objective is to extend the techniques given in \cite{MSS} for finite Abelian groups which, by the fundamental theorem of finite Abelian groups, can be written as \[
G \approx \prod_{i=1}^{N}\frac{\ZZ}{p_i^{e_i}\ZZ}
\]
For this, we will first solve the problem for finite Abelian $p$-groups and then extend the techniques to finite Abelian groups.

\subsection{Contributions of this work}
The following are the contributions of this work:
\begin{enumerate}
    \item We have arrived at the necessary and sufficient conditions for the existence of a solution to the root extraction
    problem in finite Abelian $p$-groups.
    \item We have provided an algorithm for root extraction in finite Abelian $p$-groups. The SageMath implementation of which can be found in \url{https://github.com/uacharjee14/Root-Extraction}.
    \item We have extended the algorithm for extracting roots in finite Abelian $p$-groups to finite Abelian groups.
    \item We conclude that root extraction in finite Abelian groups is no harder than solving discrete logarithms and computing basis. 
\end{enumerate}
\section{Finite Abelian $p$-groups}

The group $G$ that we consider in this section is a finite Abelian $p$-group having the structure
\begin{equation}\label{struct:finite Abelian-p}
G \approx \prod_{i=1}^{N}\frac{\ZZ}{p^{e_i}\ZZ}
\end{equation}
The $e_i$'s need not be distinct but we will assume without loss of generality that $e_j \leq e_{j+1}$. Also, we would denote the additive identity of the group $G$ by $0$.\\
The concept of \emph{basis} of a finite Abelian group has been defined in \cite[Definition 9.1]{andrewthesis}. We state it here in additive notation.
\begin{definition}[Basis of a finite Abelian group]
    A basis for a finite Abelian group $G$ is an ordered set of group elements, $\{Q_1,\ldots, Q_N\}$ with the property that every $B \in G$ can be uniquely expressed in the form $B=b_1Q_1+\cdots+b_NQ_N$ with $0\leq b_i< |Q_i|$ for $1\leq i\leq N$.
\end{definition}
\begin{remark}
    Equivalently, it can be stated that a set $\{Q_1,\ldots,Q_N\}$ is a basis of $G$ if \begin{enumerate}
        \item The set $\{Q_1,\ldots,Q_N\}$ is a \emph{spanning} set i.e., every element $B\in G$ can be expressed as $B=b_1Q_1+\cdots+b_NQ_N$ with $0\leq b_i< |Q_i|$ for $1\leq i\leq N$.
        \item The set $\{Q_1,\ldots,Q_N\}$ is \emph{linearly independent}  i.e., if $a_1Q_1+\cdots+a_NQ_N=0$ then $a_i\equiv 0\mod |Q_i|$ for $1\leq i\leq N$.
    \end{enumerate}
\end{remark}
We formulate the \emph{Root Extraction Problem} (REP) for groups of the form (\ref{struct:finite Abelian-p}) as follows:
\begin{problem}[Root Extraction Problem in finite Abelian $p$-group] \label{REP:Abelain-p}
    Let $G$ be a finite Abelian $p$-group. Given $m_1,\ldots,m_N$ where $m_i\in \frac{\ZZ}{p^{e_i}\ZZ}$  and $K \in G$, find a basis $P_1,\ldots,P_N$ of $G$ such that $K = m_1P_1+\cdots+m_NP_N $ ($|P_i|=p^{e_i}$).
\end{problem}

We introduce some definitions and results that would be useful in solving the problem. 

Let $\{Q_1,\ldots,Q_N\}$ be a basis of $G$ (that we might have computed using the algorithm in \cite[\S4]{Sutherland2008StructureCA}). Additionally let $\{Q_1,\ldots,Q_N\}$ be sorted i.e., $|Q_1|\leq |Q_2| \leq \cdots \leq|Q_N|$.

\begin{definition}[Primitive Element]
    We call an element $K=q_1Q_1+\cdots+q_NQ_N \in G$ a \emph{primitive} element if  $p \nmid q_i$ for some $1 \leq i \leq N$.
\end{definition}
\begin{proposition}
     If $K\in G$ is primitive w.r.t. some basis $\{Q_1,\ldots,Q_N\}$ then $K$ is primitive w.r.t. any other basis of $G$.
\end{proposition}
\emph{Proof: }Assume that $K$ is not primitive relative to some basis $\{P_1,\ldots,P_N\}$ of $G$. Then $K=p(m_1'P_1+\cdots+m_N'P_N)=pK'$, now we may write $K'= q_1'Q_1+\cdots+q_N'Q_N$. Hence by the uniqueness of representation $q_i = pq_i'$ for all $i$. This contradicts our hypothesis that $K$ is primitive w.r.t. $\{Q_1,\ldots,Q_N\}$.\hfill $\square$\\
\par
The above result says that the property of an element being \emph{primitive} is independent of the basis used in its representation.\\
Next, we define a function $\nu_p:G\rightarrow \ZZ_{\geq 0}\cup\{\infty\}$ as follows:
\begin{definition}
    Let a non-zero element $K\in G$ be written as $K=q_1Q_1+\cdots+q_NQ_N$. 
    The function $\nu_p:G\rightarrow\ZZ_{\geq 0}\cup\{\infty\}$ is defined as $\nu_p(K)=\nu_p(q_1Q_1+\cdots+q_NQ_N)=r$ where $r$ is the highest non-negative integer such that $p^r\mid q_i, \forall i$. Also, by convention, $\nu_p(0)=\infty$.
\end{definition}

\begin{remark}
    [The function $\nu_p$ is well defined] For a non-zero element $K=q_1Q_1+\cdots+q_NQ_N=m_1P_1+\cdots+m_NP_N \in G$, where $\{P_1,\ldots,P_N\}$ is any other basis of $G$, $\nu_p(q_1Q_1+\cdots+q_NQ_N)=\nu_p(m_1P_1+\cdots+m_NP_N)$.

\end{remark}
\emph{Proof: } Now $K = p^rK'$, where $K'$ is primitive w.r.t $\{Q_1,\ldots,Q_N\}$, so it is primitive w.r.t $\{P_1,\ldots,P_N\}$ as well. Let $K' = m_1'P_1+\cdots+m_N'P_N$ where $p \nmid m_i'$ for some $i$. Then, $\nu_p(m_1P_1+\cdots+m_NP_N)=\nu_p(p^rm_1'P_1+\cdots+p^rm_N'P_N)=r$. \hfill $\square$\\

\begin{proposition}
Following are some properties of the function $\nu_p$.
\begin{enumerate}
    \item $\nu_p(A)=\infty$ if and only if $A=0$;
    \item $\nu_p(-A)=\nu_p(A)$;
    \item $\nu_p(A+B)\geq \min(\nu_p(A),\nu_p(B))$.
\end{enumerate}
\end{proposition}
\emph{Proof:} These properties follow directly from the definition of $\nu_p$. Since properties 1 and 2 are straightforward we only prove property 3 here. Let \[
A = a_1Q_1+\cdots+a_NQ_N \text{ and } B = b_1Q_1+\cdots+b_NQ_N.
\]  Let, $\nu_p(A)=r_A$ and $\nu_p(B)=r_B$ where $r_A\leq r_B$. So, $p^{r_A}|a_i$ and $p^{r_A}|b_i$ for $1\leq i\leq N$. This means $p^{r_A}|(a_i+b_i)$ for $1\leq i\leq N$ and therefore $\nu_p(A+B)\geq r_A = \min(\nu_p(A),\nu_p(B))$.\hfill$\square$

From this, we can see that the function defined above satisfies the properties of a group valuation function defined in \cite[Chapter 2, \S 2.2]{aschenbrenner2019asymptotic}. 

Given an element $K \in G$ we say that we can \emph{extend} it to a basis of $G$ if we can find elements $\{K_2,\ldots,K_N\}$ such that $\{K,K_2,\ldots,K_N\}$ is a basis for $G$. The process of extending an element $K$ to a basis of $G$ will be referred to as \emph{basis extension}.

\begin{theorem}[Necessary Condition for Basis Extension]
     If an element $Q=q_1Q_1+\cdots+q_NQ_N \in G$ can be extended to a basis of $G$ then $Q$ is primitive.
\end{theorem}
\emph{Proof: }
We show that if $Q$ is not \emph{primitive} then it cannot be extended to a basis of $G$. Assume to the contrary that $\{Q,\QQ_2,\ldots,\QQ_N\}$ is a basis for $G$. Then since $p\mid q_k \forall k$ as $Q$ is not primitive,  we have $Q=pQ'$. Let $Q'=q_1'Q+q_2'\QQ_2+\ldots+q_N'\QQ_N$. So we have, \[
    Q=pq_1'Q+pq_2'Q_2+\cdots+pq_N'Q_N\\
    \implies (pq_1'-1)Q+pq_2'Q_2+\cdots+pq_N'Q_N=0. \]
    But this would mean $pq_1'-1 = 0$ (from linear independence) which in turn means $p$ is a unit in $\ZZ_{|Q|}$. However, $|Q|$ is a power of $p$ so $p$ cannot be a unit in $\ZZ_{|Q|}$ and therefore this is a contradiction.\hfill $\square$
\par
\begin{theorem}[Sufficient Condition for Basis Extension]\label{result:basis_extension}
    An element $Q=q_1Q_1+\cdots+q_NQ_N \in G$ can be extended to a basis of $G$ if $Q$ is primitive and  $|Q|= p^{e_k}$ where $k$ is the largest index such that $p \nmid q_k$.
\end{theorem}
\emph{Proof: }Assume that the said conditions hold, i.e. $Q$ is \emph{primitive} and $|Q|=p^{e_k}$, where $k$ is the highest index such that $p\nmid q_k$. Next, we show that the set $\{Q_1,\ldots, Q_{k-1}, Q, Q_{k+1},\ldots, Q_N\}$ is a linearly independent and spanning set. Let, \begin{equation}\label{value of q}
    Q = q_1Q_1+\cdots+q_kQ_k+\cdots+q_NQ_N
\end{equation} where $p\nmid q_k$ and $|Q|=p^{e_k}$.
\begin{enumerate}
    \item \emph{Linear independence}: Let,\[
    a_1Q_1+\cdots+a_kQ+\cdots+a_NQ_N=0.
    \] Substituting from \ref{value of q} we have,\[
    (a_1+a_kq_1)Q_1+\cdots+a_kq_kQ_k+\cdots+(a_N+a_kq_N)Q_N=0.
    \] From linear independence of basis $\{Q_1,...,Q_N\}$ we must have, $p^{e_k}\mid a_kq_k$ and $p^{e_i}\mid (a_i+a_kq_i)$ for $i\neq k$. Since, $p\nmid q_k$ so $p^{e_k}\mid a_k$ i.e., $a_k\equiv0\mod|Q|$. This also means that $a_k\equiv0\mod|Q_i|$ for all $i<k$ and therefore $a_i\equiv0\mod|Q_i|$. Also, because $|Q|=p^{e_k}$ so for all $i>k$, $p^{e_i-e_k}\mid q_i$ and therefore $q_ia_k\equiv 0\mod|Q_i|$. This means $a_i\equiv0\mod|Q_i|$ for $i>k$. This establishes linear independence.
    \item \emph{Spanning}: Let $B\in G$. Since $\{Q_1,...,Q_N\}$ is a basis of $G$ so we must be able to write $B$ as, \[
    B = b_1Q_1+\cdots+b_NQ_N.
    \] Also, since $p\nmid q_k$ in \ref{value of q} so $q_k^{-1}$ exists in $\frac{\ZZ}{p^{e_i}\ZZ}$ for $1\leq i\leq N$. This means we can write $Q_k$ as,\[
    Q_k = q_k^{-1}(Q-q_1Q_1-\cdots-q_{k-1}Q_{k-1}-q_{k+1}Q_{k+1}-\cdots-q_NQ_N).
    \] This can be substituted in the expression of $B$ to get a representation in terms of $\{Q_1,\ldots, Q_{k-1}, Q, Q_{k+1},\ldots, Q_N\}$.\hfill$\square$
\end{enumerate}
\begin{remark}
    The term \emph{primitive element} of a group has been used in \cite{CliffordGoldstein+2010+601+611} to mean an element that is part of some basis of a group. A primitive element $K \in G$ where $G$ is an Abelian $p$-group of the form (\ref{Struct2}) then $|K|=p^e$, so by theorem \ref{result:basis_extension} $K$ must be a part of some basis. 
    Therefore, our definition of primitive element for Abelian $p$-groups of the form (\ref{Struct2})  is equivalent to that in \cite{CliffordGoldstein+2010+601+611}. However, this is not the case for an arbitrary finite Abelian $p$-group, where an element of the group has to satisfy the sufficient conditions mentioned above to be a part of a basis.
\end{remark} 

\begin{example}[A primitive element that is not a part of some basis]
    Let $G\approx \frac{\ZZ}{2\ZZ}\times\frac{\ZZ}{4\ZZ}\times\frac{\ZZ}{8\ZZ}$ and let the primitive element be $Q = (1,0,2)$. Notice that $|Q|=4$ and it does not satisfy the condition in theorem \ref{result:basis_extension}. 
    Assume to the contrary that $\{Q_1,Q,Q_3\}$ is  a basis of $G$. Now,\[
    2Q=(0,0,4)=4\cdot(0,0,1)
    \]
    Let, $(0,0,1)=q_1Q_1+q_2Q+q_3Q_3$. So,\[
    2Q= 4\cdot(0,0,1) = 4\cdot (q_1Q_1+q_2Q+q_3Q_3) = 4\cdot q_1Q_1+4\cdot q_3Q_3
    \]
    Note that at least one of $4\cdot q_1$ and $4 \cdot q_3$ is non-zero.\\
    This goes against the linear independence of $\{Q_1, Q, Q_3\}$ and therefore it cannot be a basis for $G$. This is a contradiction.\hfill $\square$
\end{example}

We will now define what we mean by \emph{reducibility} and \emph{reduced form} of a root extraction problem:
\begin{definition}\label{Def:Reduced Form}[Reducibility and Reduced Form]Let, $K=q_1Q_1+\cdots+q_NQ_N \in G$. If the REP has the parameters, $K$ and $m_1,\ldots,m_N$ where $m_i \in \frac{\ZZ}{p^{e_i}\ZZ}$ then the problem is said to be \emph{reducible} if $\nu_p(K)=\nu_p(\sum_{i=1}^N m_iQ_i)$. Further, the new REP with parameters $K' \in G$ and $m_1',\ldots,m_N'$, where $K'=\sum_{i=1}^Nq_i'Q_i, q_i=p^{\nu_p(K)}q_i'$ and $m_i=p^{\nu_p(K)}m_i'$ is said to be in \emph{reduced form}. Additionally in the reduced form, we require that $q_i'=0$ when $q_i=0$ and $m_i'=0$ when $m_i=0$.
\end{definition}
 In our algorithms, we will only solve the reduced form of the REP so we state the following result:
 \begin{theorem}[Reduced form can be solved $\iff$ REP can be solved]
    Let $G$ be a finite Abelian $p$-group with a sorted basis basis $\{Q_1,\ldots,Q_N\}$ and the root extraction problem be as defined in problem \ref{REP:Abelain-p} where $\nu_p(K)=r$ for some $r \geq 0$. A solution for REP exists if and only if a solution for the reduced form exists.
\end{theorem}
\emph{Proof: }\begin{enumerate}
    \item\textbf{Reduced form can be solved  $\implies$ REP can be solved: } Suppose we reduce $(K,m_1,\ldots,m_N)$ to $(K',m_1',\ldots,m_N')$. Let $\{P_1,\ldots,P_N\}$ be a solution to the reduced from of the REP. So, $K'=m_1'P_1+\cdots+m_N'P_N$. Now, $K=p^rK'=p^rm_1'P_1+\cdots+p^rm_N'P_N=m_1P_1+\cdots+m_NP_N$. The solution is $\{P_1,\ldots, P_N\}$ in both cases.
    \item\textbf{REP can be solved$\implies$ Reduced form can be solved: } Suppose $(K,m_1,\ldots,m_N)$ is reduced to $(K',m_1',\ldots,m_N')$. Let $\{P_1,\ldots,P_N\}$ be a solution to the REP. So, $K=\sum_{i=1}^Nq_iQ_i=\sum_{i=1}^Nm_iP_i$. Also, $K'=q_1'Q_1+\cdots+q_N'Q_N$ from the reduced REP and let $K''=m_1'P_1+\cdots+m_N'P_N$ then clearly, $p^rK'=p^rK''=K$, (as $r=\nu_p(K)$), which implies $p^r(K'-K'')=0$. Let $K'-K''=R$, so $|R|\leq p^r$. By definition,  $K'=K''+R=m_1'P_1+\cdots+m_k'P_k+\cdots+m_N'P_N+R$, where $k$ is the maximum index such that $p\nmid m_k'$, so $m_k'^{-1}$ exists in $\frac{\ZZ}{|R|\ZZ}$. We also have $|P_k|=p^{e_k}>p^r$ because $m_k' \neq 0$ and therefore $p^rm_k'P_k=m_kP_k\neq0$. So from the sufficient condition of basis extension (theorem \ref{result:basis_extension}) $P_k+m_k'^{-1}R$ can be extended to a basis of $G$ which is, $\{P_1,\ldots,P_{k-1},P_k+m_k'^{-1}R,P_{k+1},\ldots,P_N \}$. Therefore we have $K' = m_1'P_1+\cdots+m_k'(P_k+m_k'^{-1}R)+\cdots+m_N'P_N$. Here, the solution of the \emph{reduced form} need not be the same as that of REP.\hfill $\square$
\end{enumerate}
Each time we reduce the problem we would need to find the highest index such that $p\nmid q_k$ and to solve the root extraction problem from there we would require the existence of $m_k^{-1}$. So we would at least need that $p\nmid m_j$ for some $j$ such that $|Q_j|=|Q_k|$. This would make sure that we can shuffle the basis $\{Q_1,\ldots,Q_N\}$ so that $m_k$ (found using the above technique) has an inverse. For this, we have the following result:
\begin{theorem}\label{result:shuffle}
    Let $K_1=q_1Q_1+\cdots+q_NQ_N$ and $K_2=m_1P_1+\cdots+m_NP_N$ are two \emph{primitive elements} of $G$ where $|K_1|=p^e$. Also, let $k$ be the highest index such that $p\nmid q_k$ and $l$ be the highest index such that $p\nmid m_l$. If $\nu_p(p^jK_1)=\nu_p(p^jK_2), \forall j \text{ such that }  0\leq j \leq e$ then $|Q_k|=|P_l|$.
\end{theorem}
\emph{Proof: } If $|Q_k|\neq |P_l|$, then without loss of generality we may assume that $|P_l|>|Q_k|$. Let, $|P_l|=p^{e_l}$ and $|Q_k|=p^{e_k}$, so $e_l > e_k$. Now, $p^{e_l-1}K_1=p^{e_l-1}(q_1Q_1+\cdots+q_NQ_N)$, and $p^{e_l-1}K_2=p^{e_l-1}(m_1P_1+\cdots+m_NP_N)$. But then we have $\nu_p(p^{e_l-1}K_1)\neq \nu_p(p^{e_l-1}K_2)$, because $p\mid q_i, \forall i>k$, whereas $p \nmid m_l$. This contradicts the condition that we started with. \hfill $\square$

\subsection{Solving for a Simpler Structure}
We will first solve the problem for the structure 
\begin{equation}
G \approx \frac{\ZZ}{p^{e_1}\mathbb{Z}}\times \frac{\ZZ}{p^{e_2}\ZZ}    
\end{equation}

where $e_1\leq e_2$.
The root extraction problem may be stated as:
\begin{problem}
    Given $G$ as mentioned above, $K \in G$, $m_1 \in \frac{\ZZ}{p^{e_1}\ZZ}$ and $m_2 \in \frac{\ZZ}{p^{e_2}\ZZ}$ find a basis $P_1,P_2$ such that $K = m_1P_1 + m_2P_2$.
\end{problem}

We will state the necessary and sufficient conditions for the existence of a solution of the REP for this structure. We will call this the existence theorem.
\begin{theorem}[Existence Theorem] The solution to the REP problem exists if and only if 
\begin{enumerate}
    \item The REP is reducible (See definition \ref{Def:Reduced Form}),
    \item $|K|=|m_1Q_1+m_2Q_2|$ for some basis $Q_1,Q_2$ of $G$.
\end{enumerate}
\end{theorem}
\emph{Proof: } Let $\{P_1,P_2\}$ be the solution of the REP. Then, $K=m_1P_1+m_2P_2$ and $\nu_p(m_1P_1+m_2P_2)=\nu_p(m_1Q_1+m_2Q_2)$ because the value of $\nu_p$ only depends on the coefficients $m_1,m_2$. So, the REP is reducible. Also, $|K|=|m_1P_1+m_2P_2|=|m_1Q_1+m_2Q_2|$ follows from the definition of order of an element and the fact that $\{Q_1,Q_2\}$ and $\{P_1,P_2\}$ are sorted bases for the group $G$. 
Now, we give a way to construct a solution when the two conditions are met. Assume that the two conditions of the \emph{Existence Theorem} hold and we have reduced our inputs to $K=q_1Q_1+q_2Q_2,m_1,m_2$. Note that we have $|K|=|m_1Q_1+m_2Q_2|$ from the conditions 1 and 2, so if $k$ is the highest index such that $p\nmid q_k$ then $k$ is the highest index such that $p\nmid m_k$ (when $e_1=e_2$ the basis $\{Q_1,Q_2\}$ may be shuffled to get this). This would mean that $m_k^{-1}$ exists in $\frac{\ZZ}{|Q_k|\ZZ}$.  We will consider different cases and solve each of them.
\begin{description}
    \item[Case 1:] $|K|=p^{e_2}$ (equivalently, $p \nmid m_2$).  Then the solution is $\{Q_1,m_2^{-1}(K-m_1Q_1)\}$.
    \item[Case 2:] $|K|=|m_1Q_1+m_2Q_2|=p^{e_1}$ (equivalently,   $p^{e_2-e_1}\mid m_2$ and $p \nmid m_1$). Then the solution is $\{m_1^{-1}(K-m_2Q_2), Q_2 \}$;
    \item[Case 3:] $|K|=p^e, e_1<e<e_2$ (equivalently, $p^{e_2-e}\mid m_2, p^{e_2-e+1}\nmid m_2, p\nmid m_1$). Then the solution is $\{m_1^{-1}q_1Q_1, (m_2')^{-1}q_2'Q_2\}$ where $m_2=p^{e_2-e}m_2', q_2=p^{e_2-e}q_2'$.\hfill$\square$
\end{description} 

\begin{remark}
    When $e_1=e_2$ this way of constructing the solution is the same as the one given in \cite{MSS}.

We will extend this technique to an arbitrary finite Abelian $p$-group. For this:
\begin{enumerate}
    \item A slightly stronger version of the \emph{existence} conditions mentioned is required.
    \item The number of cases considered in this way will be too large and this method will turn out to be clumsy. We need to bring down the number of cases considered by grouping the cases that can be handled in a similar way.
\end{enumerate}
\end{remark} 
\subsection{Solving the REP for finite Abelian $p$-groups}
We will consider the structure given in \ref{struct:finite Abelian-p}. Let us first state and prove the \emph{existence theorem} for this structure. The proof we give is constructive, and our algorithm to extract roots will readily follow.
We state a lemma before the existence theorem.
\begin{lemma}\label{existence lemma}
    Let $G$ be a finite Abelian $p$-group with basis $\{Q_1,\ldots,Q_N\}$. Also, let $K=\sum_{i=1}^Nq_iQ_i$ and $K'=\sum_{i=1}^Nq_i'Q_i$ be two elements of $G$ such that $K=p^eK'$. If $q_i= 0 \implies q_i'=0$ then $\nu_p(K)=e+\nu_p(K')$.
\end{lemma}
\emph{Proof:} Now, $p^{\nu_p(K')}\mid q_i'$ for all $1\leq i \leq N$ and $q_i=p^eq_i'$ from the theorem statement. So, $p^{e+\nu_p(K')}\mid q_i$ for all $1\leq i \leq N$. This means, $\nu_p(K)\geq e+\nu_p(K')$.\\
Note that since $K=p^eK'$ so $\nu_p(K)\geq e$. Now, $p^{\nu_p(K)}\mid q_i$ i.e., $p^{\nu_p(K)}\mid p^eq_i'$ for all $1\leq i \leq N$. Since, we have $q_i= 0 \implies q_i'=0$ i.e., $q_i\neq 0 \iff q_i'\neq0$ so $p^{\nu_p(K)-e}\mid q_i'$  for all $1\leq i \leq N$. Therefore $\nu_p(K')\geq \nu_p(K)-e$ i.e., $\nu_p(K)\leq e+\nu_p(K')$.\hfill$\square$
\begin{remark}
 Let $K=\sum_{i=1}^Nq_iQ_i\in G$ and $K=p^eK'$ where $K'=\sum_{i=1}^Nq_i'Q_i$ then it is not in general true that \begin{equation}
     \nu_p(K)=e+\nu_p(K')  \label{valuation prop}
 \end{equation} One example could be $K=(0,4) \in \frac{\ZZ}{2\ZZ}\times\frac{\ZZ}{8\ZZ}$ and $K'=(1,2)$ such that $K=2K'$.     
\end{remark}
 
\begin{theorem}[Existence Theorem]
  Let $G$ be a finite Abelian $p$-group with the structure given in \ref{struct:finite Abelian-p} with basis $\{Q_1,\ldots,Q_N\}$. Also let the root extraction problem (REP) be as defined in problem \ref{REP:Abelain-p}. Suppose $K$ can be written as $K=q_1Q_1+\cdots+q_NQ_N$ and $|K|=p^e$. Then the solution to REP exists if and only if
\begin{enumerate}
    \item $|K|=|(m_1Q_1+\cdots+m_NQ_N)|$, and 
    \item $\nu_p(p^jK)=\nu_p(p^j(m_1Q_1+\cdots+m_NQ_N)), \forall j \text{ such that }0\leq j<e$.
\end{enumerate}  
\end{theorem}
\emph{Proof: }
\begin{description}
    \item[If a solution to REP exists, then the properties hold. ]If $K=\sum_{i=1}^Nm_iP_i$, then $|K|=|\sum_{i=1}^Nm_iP_i|=|\sum_{i=1}^Nm_iQ_i|$ and also $\nu_p(p^jK)=\nu_p(p^j(m_1P_1+\cdots+m_NP_N))=\nu_p(p^j(m_1Q_1+\cdots+m_NQ_N))$ for $0\leq j<e$.
    \item[If the properties hold a solution to REP exists. ] 
We will use induction on the rank of $G$ calling it $n$. Verifying that a solution exists for $n=1$ if the properties hold is easy. For $n=2$ we already constructed a solution in the previous subsection. Also, it is to be noted that for all the cases that we have solved the problem whenever $m_i=0$ we had $P_i=Q_i$ i.e., the $i$-th \emph{basis} element was not altered. We now assume that for all $n<N$ the problem can be solved and try to construct a solution for $n=N$.\\
As $\nu_p(K)=\nu_p((m_1Q_1+\cdots+m_NQ_N))$ so the problem is \emph{reducible}. Let $K', m_1',\ldots,m_N'$ be the \emph{reduced problem}. If $K'=\sum_{i=1}^Nq_i'Q_i$, let $k$ be the greatest index such that $p\nmid q_k'$ and so $p \nmid m_k'$ (assuming the necessary shuffling of basis is done using result \ref{result:shuffle}). Construct a set of indices $I_Q$ that contains all indices $i$ for which $p^{e_k}q_i'Q_i=0$. Similarly, construct a set of all indices $I_M$ that contains all indices for which $p^{e_k}m_i'Q_i=0$.\\
\par
Define a new subproblem with parameters $K_1=\sum_{i\in I_Q}q_iQ_i$ and $\mm_i$ where $\mm_i=m_i, i\in I_M$, otherwise $\mm_i=0$. Notice that this problem is reducible as well because $\nu_p(K)= \nu_p(K_1)$ and $\nu_p((m_1Q_1+\cdots+m_NQ_N))=\nu_p(\sum_{i=1}^{N}\mm_iQ_i)$. We reduce this problem to obtain the parameters $K_1'=\sum_{i\in I_Q}q_i'Q_i$ and $\mm_i'$. This subproblem has a straightforward solution where $P_i=Q_i, i \neq k$, and $P_k= \mm_k'^{-1}(K_1'-\sum_{i \neq k}\mm_i'Q_i)$ . This solves the subproblem and we obtain our new \emph{basis} $\{P_1,\ldots,P_N\}$ (this follows from result \ref{result:basis_extension}) such that $K_1=\sum_{i=1}^{N}\mm_iP_i=\sum_{i\in I_M}m_i'P_i$. If $e_k=e$ then we are done and there is no need to define another subproblem.\\
\par
When $e_k\neq e$ the other subproblem can be defined with the parameters $K_2=K-K_1=\sum_{i\notin I_Q}q_iQ_i$ and $\mm_i=m_i, i\notin I_M$, otherwise $\mm_i=0$. Notice that for all $i\leq k$ we have $\mm_i=0$ and similarly, $K_2$ can be regarded as $K_2=\sum_{i=1}^N\qq_iQ_i$ where $\qq_i=q_i, i\notin I_Q$, otherwise $\qq_i=0$. This problem is equivalent to solving a problem where $n=N-k$ and our mentioned conditions hold. 
\begin{enumerate}
    \item As, $|K|=|K_2|$ and $|(m_1Q_1+\cdots+m_NQ_N)|=|(\mm_{k+1}Q_{k+1}+\cdots+\mm_{N}Q_N)|$. So, $|K_2|=|(\mm_{k+1}Q_{k+1}+\cdots+\mm_NQ_N)|$
    \item Because $\nu_p(p^jK)=\nu_p(p^j(m_1Q_1+\cdots+m_NQ_N))$ for $0\leq j<e$ so we have from lemma \ref{existence lemma} $\nu_p(p^jK_2)=\nu_p(p^j(\mm_{k+1}Q_{k+1}+\cdots+\mm_NQ_N))$ for $1\leq p^j < |K_2|$.
\end{enumerate}
Therefore, from the induction hypothesis we have a solution for this problem as well say $\{\PP_{k+1},\ldots,\PP_N\}$ with the additional condition that $\PP_i=P_i$ whenever $\mm_i=0$. So, taking $\PP_i=P_i, \forall i\leq k$ we have the \emph{basis} set $\{\PP_{1},\ldots,\PP_N\}$. Now , $\sum_{i=1}^Nm_i\PP_i=\sum_{i\in I_M}m_i\PP_i+\sum_{i\notin I_M}m_i\PP_i=K_1+K_2=K$.\hfill $\square$.
\end{description}

\subsection{The Algorithm for root extraction in finite Abelian $p$-groups}
We have divided the algorithm into three parts:
\begin{description}
    \item[Algorithm \ref{alg: check}: ] An algorithm that returns a solution if the existence conditions are satisfied.
    \item[Algorithm \ref{alg:Reduce}: ] An algorithm to convert the given inputs to \emph{reduced form}.
    \item[Algorithm \ref{alg: RootExtract}: ] The algorithm for root extraction.
\end{description}
We have assumed that the following are available globally and can be accessed and modified by each algorithm that we have mentioned :
\begin{enumerate}
    \item The group description $G$,
    \item $Q_1,\ldots,Q_N$ (the computed \emph{basis} of group $G$), 
    \item $q_1,q_2,\ldots,q_N$ ( the discrete logarithm of $K$ w.r.t $\{Q_1,\ldots,Q_N\}$),
    \item $m_1,\ldots,m_N$ (the coefficients supplied as inputs to the root extraction problem).
\end{enumerate}
\hfill\break
\begin{algorithm}[!ht]
\caption{Checking if Root Extraction is possible}
    \label{alg: check}
Input:$G, K, m_1,\ldots,m_N$\\
Precondition:$K\in G, m_i \in \mathbb{Z}/p^{e_i}\mathbb{Z}, e_1\leq e_2\leq\ldots\leq e_N$\\
Output:$P_1,\ldots,P_N$\\
Postcondition:$K=m_1P_1+\cdots+m_NP_N$ and $\langle P_1,\ldots,P_N \rangle= G$\\
\textbf{begin}
\begin{algorithmic}[1]
    \IF{$m_1=\ldots=m_N=0 \OR q_1=\ldots=q_N=0$}
        \IF{$(\exists : q_i \neq 0)\OR (\exists i :m_i \neq 0)$}
            \STATE Raise Exception(Solution doesn't Exist) 
        \ELSE
            \STATE return $Q_1,\ldots,Q_N$
        \ENDIF
    \ENDIF
    \IF{$|K|=|m_1Q_1+\cdots+m_NQ_N| \AND\nu_p(p^{e}q_1Q_1+\cdots+p^{e}q_NQ_N)=\nu_p(p^{e}m_1Q_1+\cdots+p^{e}m_NQ_N) ,0\leq e < \log_p(|K|)$}
        \STATE Pass parameter $K$ to Algorithm \ref{alg: RootExtract}  
        \RETURN $P_1,\ldots,P_N$ 
    \ELSE 
        \STATE raise Exception\{Solution doesn't exist\}
        \RETURN 
    \ENDIF
\end{algorithmic}
\end{algorithm}

\begin{algorithm}
    \caption{Reduction algorithm:}
    \label{alg:Reduce}
    Input: $K, I_Q, I_M$\\
    Precondition: $K\in G, m_i \in \mathbb{Z}/p^{e_i}\mathbb{Z}, e_1\leq e_2\leq\ldots\leq e_N$ (The conditions from Algorithm \ref{alg: check} are satisfied).\\
    Output: $K$, (Also the variables $q_1,\ldots,q_N$ and $m_1,\ldots, m_N$ are altered)\\
    Post-condition: The problem is reduced (i.e. $\nu_p(K)=\nu_p(m_1Q_1+\cdots+m_NQ_N)=0$)\\
    \textbf{ Reduce($K,I_Q,I_M$):}
    \begin{algorithmic}[1]
        \STATE find largest $r$ such that $p^r\mid q_i ,\forall i \in I_Q$ 
        \IF{$r\neq0$}
            \STATE set $q_i = q_i/p^r, \forall i \in I_Q$ 
            \STATE set $m_i = m_i/p^r, \forall i \in I_M$ 
            \STATE $K=\sum_{i\in I_Q}(q_i/p^r)Q_i$ 
        \ENDIF
        \RETURN $K$
    \end{algorithmic}
\end{algorithm}

\begin{algorithm}
    \caption{Root Extraction Algorithm for finite Abelian p-groups}
    \label{alg: RootExtract}
    Input: $K$\\
    Precondition: Since the conditions have already been checked by Algorithm \ref{alg: check}, we do not specify any precondition here\\
    Output: $P_1,\ldots,P_N$\\
    Post-condition: $K=m_1P_1+\cdots+m_NP_N$ and $\langle P_1,\ldots,P_N \rangle= G$ (Same as that in Algorithm \ref{alg: check})\\
    \textbf{begin}
    \begin{algorithmic}[1]
        \STATE set $I_M=I_Q=\{1,\ldots,N\}$
        \STATE set $P_i=Q_i \forall i \in \{1,\ldots,N\}$
        \WHILE{$\TRUE$}
            \STATE $K=Reduce(K,I_Q,I_M)$ 
            \STATE Find maximum index $k$ such that $p \nmid m_k, k\in I_M$ 
            \IF{$|K|=p^{e_k}$}
                \STATE set $P_k = m_k^{-1}(K-\sum_{i \in I_M, i \neq k}m_iQ_i)$ 
                \RETURN
            \ELSE
                \STATE set $I_M'=\{j:p^{e_k}m_jQ_j=0\, j \in I_M\}$
                \STATE set $I_Q' = \{
                j:p^{e_k}q_jQ_j=0, j \in I_Q\}$ 
                \STATE set $K'= \sum_{i \in I_Q'}q_iQ_i$ 
                \STATE set $P_k= m_k^{-1}(K'-\sum_{i \in I_M', i \neq k}m_iQ_i)$ 
                \STATE update $I_M = I_M-I_M'$ 
                \STATE update $I_Q = I_Q-I_Q'$ 
                \STATE update $K = K-K'$
            \ENDIF
        \ENDWHILE
    \end{algorithmic}
\end{algorithm}
\textbf{Number of group operations required in Algorithm \ref{alg:Reduce}}\\
Only in step 5, we are performing group operations. We are exponentiating at most $N$ times and adding at most $N$ times. This is bounded by $O(Ne_N\log_2p+N)$ group operations which in turn is bounded by $O(Ne_N\log_2p)$.\\
\par
\textbf{Number of group operations required in Algorithm \ref{alg: RootExtract}}\\
We first find a bound on the number of group operations in one iteration. In step 4 we call Algorithm \ref{alg:Reduce} and it is bounded by $O(Ne_N\log_2p)$ group operations. After having found 
 the index $k$  in step 5, to check if the order of $K$ is $p^{e_k}$ in step 6 requires one exponentiation which is bounded by $O(e_N\log_2p)$ group operations. Inside the if block, step 7 requires at most $N$ exponentiations and additions so it is bounded by $O(Ne_N\log_2p+N)$ operations. Inside the else block, steps 10, 11, 12 13  will have at most $N$  exponentiation( and additions in step 13) so it is bounded by  and$O(Ne_N\log_2p)$. The maximum of all these is $O(Ne_N\log_2p)$. The loop has at most $N$ iterations so the total number of group operations is bounded by $O(N^2e_N\log_2p)$ group operations.
\par
\textbf{Number of group operations required in Algorithm \ref{alg: check}}\\
Computing order of $K$ takes $O(e_N^2\log_2p)$. In step 8 $N$ additions and $N$ exponentiations take $O(Ne_N\log_2p+N)$. We need to do the same additions and exponentiations again $\log_p(|K|)$ times which is bounded by $e_N$ so the total number of group operations here is bounded by $O(Ne_N^2\log_2p)$. Algorithm \ref{alg: RootExtract} takes $O(N^2e_N\log_2p)$ group operations. So the total number of group operations required is $O((N+e_N)Ne_N\log_2p)$.\\

We have provided some examples below:
\begin{example}\label{Ext_Ex1}
    Consider the group $G \approx  \frac{\ZZ}{2\ZZ}\times \frac{\ZZ}{8\ZZ}\times \frac{\ZZ}{16\ZZ} $, with standard basis. $K = (0,2,8)$ and $m_1=0, m_2=6,m_3=4$. \\ Note that $|(0,2,8)|=|(0,6,4)|=4$, $\nu_2((0,2,8))=\nu_2((0,6,4))=1$ and $\nu_2((0,4,0))=\nu_2((0,4,8))=2$, so our conditions are satisfied. Next, to solve the problem we reduce it to:\\
    $K'=(0,1,4), m_1=0, m_2=3, m_3=2$.\\
    Now, the highest index $k$ such that $p\nmid m_k$ is 2 and note that $|K'|=2^{e_2}=2^3=8$ so the solution is:
    \\$P_1=(1,0,0), P_2=3^{-1}((0,1,4)-2(0,0,1))=(0,3,6), P_3=(0,0,1)$. \hfill $\square$
\end{example}
    
\begin{example}
    We will use the same group as in Example \ref{Ext_Ex1} with standard basis. Let $K=(1,2,2)$ and $m_1=1,m_2=6,m_3=10$. Notice that $|(1,2,2)|=|(1,6,10)|=8$, $\nu_2((1,2,2))=\nu_2((1,6,10))=0, \nu_2((0,4,4))=\nu_2((0,4,4))=2$ and $\nu_2((0,0,8))=\nu_2((0,0,8))=3$. This problem is already \emph{reduced}:
    Now, in $K$,  $2 \nmid 1$ and $|K|=8$, so in this case we will need to partition the problem. Clearly, we may write $K=K_1+K_2$ where $K_1=(1,0,0)$ and $K_2=(0,2,2)$. We must also partition our coefficients $\{1,6,10\}=\{1\}\cup\{6,10\}$. \begin{enumerate}
        \item Solving for $K_1, m_1$, i.e. $(1,0,0)$ and $m_1 = 1$. We obtain $P_1 = (1,0,0)$. Note that $P_2, P_3$ are still the standard basis.
        \item Solving the other sub-problem now i.e. $K_2, m_2, m_3$. We first reduce it to get $(0,1,1)$ and $m_2'= 3, m_3'=5$. This is can be easily solved as $P_2=(0,1,0), P_3=5^{-1}((0,1,1)-3\cdot(0,1,0))=(0,6,13)$. \hfill $\square$
    \end{enumerate}
\end{example}

\begin{example}
    Consider $G \approx \frac{\ZZ}{4\ZZ}\times\frac{\ZZ}{16\ZZ}\times\frac{\ZZ}{32\ZZ}\times\frac{\ZZ}{64\ZZ}$ with standard basis, $K=(3,2,8,4)$, and $M=[1,6,4,12]$ denote the multipliers supplied as inputs to the root extraction problem.\\
First, we ensure that our conditions are satisfied.
\begin{enumerate}
    \item Note that $|K|=|\sum_{i=1}^4m_iQ_i|=16$
    \item We need to check the $\nu_2()$ value for $p=1,2,4,8$, i.e., check that $\nu_p(pK)=\nu_p(p\sum_{i=1}^4m_iQ_i)$
        \begin{enumerate}
            \item $\nu_2(3,2,8,4)=\nu_2(1,6,4,12)=0$
            \item $\nu_2(2,4,16,8)=\nu_2(2,12,8,24)=1$
            $\nu_2(0,8,0,16)=\nu_2(0,8,16,48)=3$
            \item $\nu_2(0,0,0,32)=\nu_2(0,0,0,32)=5$
        \end{enumerate}
    \item The initial problem can be split as $K_1=(3,0,8,0), M_1=(1,0,0,0)$ and $K_2=(0,2,0,4),M_2=(0,6,4,12)$. Of these, the first one can be solved readily. We get $P_1=(3,0,8,0)$.
    \item We may further reduce and partition the second sub-problem. $K_2'=(0,1,0,2),M_2'=(0,3,2,6)$ to $K_3=(0,1,0,0), M_3=(0,3,2,0)$ and $K_4=(0,0,0,2), M_4=(0,0,0,6)$. $K_3, M_3$ can be readily solved. We obtain $P_2=3^{-1}((0,1,0,0)-2Q_3)=(0,11,26,0), P_3=Q_3=(0,0,1,0)$.
    \item $K_4, M_4$ can be reduced to $K_4'=(0,0,0,1), M_4'=(0,0,0,3)$ and solved. So $P_4=3^{-1}(0,0,0,1)=(0,0,0,43)$. \hfill $\square$
\end{enumerate}
\end{example}

We also provide examples of cases when a solution to the root extraction problem doesn't exist.
\begin{example}
    Let $G\approx \frac{\ZZ}{4\ZZ}\times\frac{\ZZ}{16\ZZ}$ with standard basis, $K=(1,0)$ and $M=(m_1,m_2)=(0,4)$.
    Now we check the existence conditions. Note that $|K|=|M|=4$. However, $\nu_2(K)=0$ whereas $\nu_2(M)=2$ so a solution to this problem does not exist. \hfill $\square$
\end{example}

\begin{example}
    Consider $G\approx \frac{\ZZ}{2\ZZ}\times\frac{\ZZ}{4\ZZ}\times\frac{\ZZ}{32\ZZ}\times\frac{\ZZ}{64\ZZ}$ with standard basis. Let $K=(0,1,2,0)$ and $M=(m_1,m_2,m_3,m_4)=(1,1,4,4)$. We check the existence conditions now. Note that $|K|=|M|=16$, $\nu_2(K)=\nu_2(M)=0$ and $\nu_2(2K)=\nu_2(2M)=1$. However, $\nu_2(4K)=\nu(0,0,8,0)=3$ and $\nu_2(4M)=\nu_2(0,0,16,16)=4$, so $\nu_2(4K)\neq \nu_2(4M)$. Therefore a solution to this problem doesn't exist.\hfill $\square$
\end{example}

 The algorithm for root extraction in finite Abelian $p$-groups [Algorithm \ref{alg: check}, \ref{alg:Reduce} and \ref{alg: RootExtract}] has been implemented in SageMath and the code can be found in \url{https://github.com/uacharjee14/Root-Extraction}.
\par
\section{Finite Abelian Groups}

We will extend our results from the previous section to solve the \emph{root extraction} problem in generic finite Abelian groups.\\
From the fundamental theorem of finite Abelian groups, the structure of a finite Abelian group $G$ is \[
G \approx \prod_{i=1}^{N}\frac{\ZZ}{p_i^{e_i}\ZZ}.
\]
We will assume that this structure is already known. This can be found using Sutherland's structure computation algorithm given in \cite[\S5]{Sutherland2008StructureCA}.
This structure can be restated as:
\[
G \approx \prod_{j=1}^{N_1}\frac{\ZZ}{p_1^{e_{1j}}\ZZ}\times \cdots\times\prod_{j=1}^{N_r}\frac{\ZZ}{p_r^{e_{rj}}\ZZ}\]
Here $e_{ij}\leq e_{ik}$ for $j \leq k$ and $p_i\neq p_j$ for $i\neq j$.
The \emph{Root Extraction Problem} can be stated as:
\begin{problem}[Root Extraction Problem in finite Abelian groups]\label{REP:finite Abelian}
    Given $m_{11},\ldots,m_{rN_r}$ such that $m_{ij} \in \frac{\ZZ}{p_i^{e_j}\ZZ}$ and $K \in G$ find a basis $\{P_{11},P_{12},\ldots,P_{rN_r}\}$ of $G$ such that $K= m_{11}P_{11}+\cdots+m_{rN_r}P_{rN_r}$.
\end{problem}
By using Sutherland's algorithm  \cite[\S5]{Sutherland2008StructureCA}) we may find a basis $\{Q_1,\ldots,Q_N\}$ for the group $G$ and then compute their orders. Then using the basis, (with known orders) we can partition the given group into Sylow $p$-subgroups and solve our problem in each of them using algorithm \ref{alg: check} discussed in the previous section.
\par
Let $\{Q_{1},Q_2,\ldots,Q_N\}$ be a basis for $G$ which can be renumbered as $\{Q_{11}, Q_{12},\ldots,Q_{rN_r} \}$ so that $|Q_{ij}|= p_i^{e_{ij}}$.\\
Then we use the Generalized Discrete Logarithm Algorithm \cite[\S3, Algorithm 2]{Sutherland2008StructureCA} to express $K$ as \[K = q_{11}Q_{11}+q_{12}Q_{12}+\cdots+q_{rN_r}Q_{rN_r}.\]
Let $G_i = \langle Q_{i1},Q_{i2},\ldots,Q_{i(N_{i})}\rangle$ and $K_i = q_{i1}Q_{i1}+\cdots+q_{i(N_{i})}Q_{i(N_{i})}$ where $1\leq i\leq r$. Now it can be seen that \[
G=G_1\oplus\cdots\oplus G_r \text{ and } K=K_1+\cdots+K_r.
\]
Note that subgroup $G_i$ is the Sylow $p_i$-subgroup of $G$ and $K_i\in G_i$ for all $i$ such that $1\leq i\leq r$.
We then pass the following inputs to Algorithm \ref{alg: check}:
\begin{enumerate} 
    \item $G_i = span\{Q_{i1},Q_{i2},\ldots,Q_{i(N_{i})}\}$ (the Group description for Sylow-$p_i$ group)
    \item $K_i = q_{i1}Q_{i1}+\cdots+q_{i(N_{i})}Q_{i(N_{i})}$ (the element of the group $G_i$)
    \item $m_{i1},\ldots,m_{iN_i}$ (the coefficients)
\end{enumerate}
This will define our $i$-th subproblem (for the prime $p_i$) \label{sub-problem}\\\\
Finally, we would make an important statement:
\begin{theorem}
    The REP for finite Abelian groups has a solution if and only if each of the sub-problems has a solution
\end{theorem}
\emph{Proof: }\begin{enumerate}
    \item \emph{If each of the sub-problems has a solution then the REP has a solution}: Let, $\beta_i=\{P_{i1},\ldots, P_{iN_i}\}$ is a solution for the $i$-th subproblem. Then, clearly $K=m_{11}P_{11}+\cdots+m_{rN_r}P_{rN_r}$. Also, $\beta=\cup_{i=1}^r\beta_i=\{P_{11},\ldots,P_{rN_r}\}$ is the \emph{basis} for the group $G$ as it is the union of the basis of all the Sylow $p_i$-subgroups of $G$.
    \item \emph{If a solution to REP exists then each of the subproblems also has a solution: } Let $\{P_{11},\ldots, P_{rN_r}\}$ be a solution to the REP. Now since $\langle P_{i1},\ldots, P_{iN_i}\rangle$ is a Sylow $p_i$-subgroup of $G$. Since, a Sylow $p_i$-subgroup of finite Abelian group is unique, so $\langle P_{i1},\ldots, P_{iN_i}\rangle=\langle Q_{i1},\ldots, Q_{iN_i}\rangle$. Therefore, $\{P_{i1},..., P_{iN_i}\}$ is a solution to the $i$-th subproblem.\hfill $\square$
\end{enumerate} 
\par
\begin{algorithm}
\caption{Root Extraction in Finite Abelian Groups}
    \label{alg: finite Abelian}
Input:$G, K, m_{11},...,m_{rN_r}$,a \emph{sorted} basis $\{Q_{11},Q_{12},...,Q_{rN_r}\}$ for $G$, $|Q_{ij}|=p_i^{e_{ij}}$, $q_{11},q_{12},...,q_{rN_r}$ such that $K= q_{11}Q_{11}+q_{12}Q_{12}+...+q_{rN_r}Q_{rN_r}$ \\
Precondition:$K\in G, m_{ij}\in \frac{\ZZ}{p_i^{e_{ij}}\ZZ}$, $e_{ij}\leq e_{ik}$ for $j \leq k$ , $p_i\neq p_j$ for $i\neq j$\\
Output:$P_{11},P_{12},...,P_{rN_r}$\\
Postcondition: $K=m_{11}P_{11}+m_{12}P_{12}+...+m_{rN_r}P_{rN_r}$ and $\{ P_{11},P_{12},...,P_{rN_r} \}$ is a basis of $G$\\
\textbf{begin}
\begin{algorithmic}[1]
    \STATE set $i=1, A=[\ ]$
    \WHILE{$i \leq r$}
        \STATE define $G' = span\{Q_{i1},Q_{i2},...,Q_{iN_i}\}$
        \STATE set $K' = q_{i1}Q_{i1}+...+q_{iN_i}Q_{iN_i}$
        \STATE pass $G',K', m_{i1},...,m_{iN_i} $ to Algorithm \ref{alg: check}.
        \STATE append output of Algorithm \ref{alg: check} to $A$
        \STATE set $i=i+1$
    \ENDWHILE
    \RETURN $A$
    
\end{algorithmic}
\end{algorithm}
\par
\textbf{The number of group operations required in Algorithm \ref{alg: finite Abelian}}\\
In step 4 we perform $N_i$ additions for each $i$. Now $N_i$ is the rank of the Sylow $p_i$-subgroup and suppose $N$ is the maximum among all such ranks i.e., $N=\max_{i=1}^rN_i$. Then, the number of group operations for additions is bounded by $O(N)$, and exponentiations take $O(e_{N_i}\log_2p_i)$ group operations. So in the loop, the total number of group operations required is $O(re_{iN_i}\log_2p_i+rN)$. Let $e=\max_{i=1}^re_{iN_i}$ and $p=\max_{i=1}^rp_i$. Algorithm \ref{alg: check} uses $O((e+N)Ne\log_2p)$. This we do for each of the $r$ primes in the structure of $G$ so the total number of group operations is bounded by $O((e+N)rNe\log_2p)$.\label{complexity}\\

\subsection{Discussion}
In all the algorithms presented till now, we have assumed that a basis $\{Q_1,\ldots, Q_N\}$ of the group and the group structure is known. We have also assumed that the discrete logarithm of the element $K$ from problem \ref{REP:finite Abelian} with respect to the basis $\{Q_1,\ldots,Q_N\}$ is known.  
The bound on the number of group operations required in the generalized discrete logarithm algorithm can be found in \cite[\S3, equation 15]{Sutherland2008StructureCA} and for the basis computation algorithm this can be found in \cite[\S5, corollary 3]{Sutherland2008StructureCA}.  
From the bound on the number of group operations required for algorithm \ref{alg: finite Abelian} we can see that \emph{root extraction in finite Abelian groups is no harder than solving discrete logarithms and computing basis}.

Further, Damgard and Koprowski \cite{cryptoeprint:2002/013} have proved that the root extraction problem in finite Abelian groups of unknown orders has an exponential lower bound. We have shown that the number of group operations for root extraction is dominated by that of basis computation and solving discrete logarithm which have exponential complexities. Thus, our results agree with that of Damgard and Koprowski. 
\bibliographystyle{acm}
\bibliography{References}
\end{document}